\documentclass{amsart}
\usepackage{amssymb}
\newcommand\R{{\mathbb R}}

\newcommand\N{{\mathbb N}}
\newcommand\mf{{\mathcal M}}
\newcommand\pad{\div}

\newtheorem{thm}{Theorem}
\newtheorem{prop}[thm]{Proposition}
\newtheorem{lem}[thm]{Lemma}
\newtheorem{cor}[thm]{Corollary}
\newtheorem{defn}[thm]{Definition}
\newcommand\emdef[1]{\emph{#1}}

\begin{document}
\title{Lower Bounds on the Growth of Grigorchuk's Torsion Group}
\author{Laurent Bartholdi}
\date\today
\email{Laurent.Bartholdi@math.unige.ch}
\address{\parbox{.4\linewidth}{Section de Math\'ematiques\\
    Universit\'e de Gen\`eve\\ CP 240, 1211 Gen\`eve 24\\
    Switzerland}}
\thanks{This work has been supported by the ``Swiss National Science Foundation''}
\begin{abstract}
  In 1980 Rostislav Grigorchuk constructed a group $G$ of intermediate
  growth, and later obtained the following estimates on its
  growth $\gamma$~\cite{grigorchuk:gdegree}:
  $$e^{\sqrt{n}}\precsim\gamma(n)\precsim e^{n^\beta},$$
  where $\beta=\log_{32}(31)\approx0.991$. He conjectured that the
  lower bound is actually tight.
  
  In this paper we improve the lower bound to
  $$e^{n^\alpha}\precsim\gamma(n),$$
  where $\alpha\approx0.5157$, with
  the aid of a computer. This disproves the conjecture that the lower
  bound be tight.
\end{abstract}
\maketitle

\section{Introduction}
The growth of finitely generated groups, in relation with properties
of differentiable manifolds and coverings, was studied since the
1950's in the former \textsc{Ussr}~\cite{svarts:growth} and in the
1960's in the West~\cite{milnor:solvable}: let the finitely generated
group $G$ be the fundamental group of a compact $CW$-complex $K$. Then
the growth of $G$ is equivalent to the growth of the universal cover
$\tilde K$.  There are well-known classes of groups of polynomial
growth: abelian ($K$ a torus, $\tilde K$ Euclidean space), and more
generally virtually nilpotent groups~\cite{gromov:nilpotent}; and
classes of exponential growth: non-virtually-nilpotent
linear~\cite{tits:linear} or non-elementary
hyperbolic~\cite{ghys-h:gromov} groups ($K$ a negatively curved
complex). John Milnor asked whether there exist finitely generated
groups with growth greater than polynomial but less than
exponential~\cite{milnor:5603}. The first example of such a group, of
intermediate growth, was discovered by Rostislav Grigorchuk,
see~\cite{grigorchuk:growth,grigorchuk:gdegree,grigorchuk:kyoto}. He
showed that the growth $\gamma$ of his group satisfies
$$e^{\sqrt{n}}\precsim\gamma(n)\precsim e^{n^\beta},$$
where
$\beta=\log_{32}(31)\approx0.991$. The author obtained the
following improvement in~\cite{bartholdi:upperbd}: Let $\eta$ be the
real root of the polynomial $X^3+X^2+X-2$, and set
$\beta'=\log(2)/\log(2/\eta)\approx0.767$. Then the growth $\gamma$ of
Grigorchuk's group satisfies
$$e^{\sqrt{n}}\precsim\gamma(n)\precsim e^{n^{\beta'}}.$$
Recently, the same result was rediscovered and extended to a larger
class of groups by Roman Muchnik and Igor Pak~\cite{muchnik-p:growth}.

A remaining outstanding problem in the theory of growth of groups is
the question, raised by Grigorchuk, of the existence of groups with
growth exactly $e^{\sqrt n}$. Such groups would have interesting
extremal properties, for instance being of finite
width~\cite{bartholdi-g:lie}. He conjectured
in~\cite{grigorchuk:gdegree} that his group has this property; but in
this paper we disprove that conjecture with the following result:
\begin{thm}\label{thm:main}
  The growth of Grigorchuk's $2$-group $G$ satisfies
  $$\gamma(n)\succsim e^{n^\alpha},$$
  where $\alpha=0.5157$.
\end{thm}
However, the following construction could produce a group of growth
$e^{\sqrt n}$: let $G$ be Grigorchuk's $2$-group, generated by
$S=\{a,b,c,d\}$. Let $\mathfrak L$ be the graded Lie $2$-algebra of $G$
(see~\cite{bartholdi-g:lie}), and let $\mathfrak G$ be the group
generated by the $1+(s-e)$, $s\in S$ in the enveloping algebra
$\mathcal U(\mathfrak L)$. There are interesting relations between the
growth of $G$ and that of $\mathfrak G$, and it seems the growth of
$\mathfrak G$ is ``smoother''.

The approach used in this paper to obtain lower bounds on the growth
was suggested, with small variations, by Yuri\u\i\ 
Leonov~\cite{leonov:pont}, where he announced $\gamma(n)\succsim
e^{n^\beta}$ for $\beta=\log_{87/22}(2)\cong0.504$. I wish to thank
Yuri\u\i for introducing me to this question, and Slava Grigorchuk and
Pierre de la Harpe for their interest.

\section{Growth of Groups}
Let $G$ be a group generated as a monoid by a finite set $S$. A
\emdef{weight} on $(G,S)$ is a function $\omega:S\to\R_+^*$. It
induces a norm (again called a \emdef{weight}) $\partial_\omega$ on $G$
by
$$\partial_\omega:\begin{cases}G\to\R_+\\
  g\mapsto\min\{\omega(s_1)+\dots+\omega(s_n)|\,s_1\cdots s_n=_G
  g\}.\end{cases}$$
The important properties of the weight
$\partial_\omega$ are its \emph{submultiplicativity}:
$\partial_\omega(gh)\le\partial_\omega(g)+\partial_\omega(h)$ and its
\emph{properness}: $|\{g\in G|\,\partial_\omega(g)\le n\}|<\infty$ for
all $n\in\N$.

A \emdef{minimal form} of $g\in G$ is a representation of $g$ as
a word of minimal weight over $S$. If a minimal form has been fixed,
it will be denoted by $\overline g$. The \emdef{growth} of $G$ with
respect to $\omega$ is then
$$\gamma_\omega:\begin{cases}\R_+\to\R_+\\ n\mapsto \#\{g\in
G|\,\partial_\omega(g)\le n\}.\end{cases}$$
Alternatively, $S$ can altogether be suppressed from the definition,
and a weight can be defined as a function $\omega:G\to\R_+$ that
has finite support generating $G$ as a monoid.

The function $\gamma:\R_+\to\R_+$ is \emdef{dominated} by
$\delta:\R_+\to\R_+$, written $\gamma\precsim\delta$, if there is a
constant $C\in\R_+$ such that $\gamma(n)\leq\delta(Cn)$ for all
$n\in\R_+$. Two functions $\gamma,\delta:\R_+\to\R_+$ are
\emdef{equivalent}, written $\gamma\sim\delta$, if
$\gamma\precsim\delta$ and $\delta\precsim\gamma$.

The following lemmata are well known:
\begin{lem}\label{lem:grequiv}
  Let $S$ and $S'$ be two finite generating sets for the group $G$,
  and let $\omega$ and $\omega'$ be weights on $(G,S)$ and $(G,S')$
  respectively. Then $\gamma_\omega\sim\gamma_{\omega'}$.
\end{lem}
\begin{proof}
  Without loss of generality we may suppose $1\notin S$, so
  $\partial_\omega(s)\neq0$ for any $s\in S$.  Let $C=\max_{s\in
    S}\partial_{\omega'}(s)/\partial_\omega(s)$. Then
  $\partial_{\omega'}(g)\leq C\partial_\omega(g)$ for all $g\in G$,
  and thus $\gamma_\omega(n)\le \gamma_{\omega'}(Cn)$, from which
  $\gamma_\omega\precsim\gamma_{\omega'}$. The opposite relation holds
  by symmetry.
\end{proof}
In particular, for fixed $S$, any weight $\partial_\omega$
is equivalent to the length-weight $|g|=\min\{n|g=s_1\dots
s_n,\,s_i\in S\}$.

The \emdef{growth type} of a finitely generated group $G$ is the
$\sim$-equivalence class containing its growth functions; it will be
denoted by $\gamma_G$.

Note that all exponential functions $b^n$ are equivalent, and
polynomial functions of different degrees are inequivalent; the same
holds for the subexponential functions $e^{n^\alpha}$. We have
$$0\precnsim n\precnsim n^2\precnsim\cdots\precnsim
e^{n^\alpha}\precnsim e^{n^\beta}\precnsim\dots\precnsim
e^n\qquad\text{ for }0<\alpha<\beta<1.$$
Note also that the ordering
$\precsim$ is not linear. Actually, Slava Grigorchuk showed in his
pioneering paper~\cite[Theorem~7.2]{grigorchuk:gdegree} that $\precsim$ admits
chains and antichains with the cardinality of the continuum.

\begin{lem}
  Let $G$ be a finitely generated group. Then $\gamma_G\precsim e^n$.
\end{lem}
\begin{proof}
  Choose for $G$ a finite generating set $S$, and define the weight
  $\omega$ by $\omega(s)=1$ for all $s\in S$. Then
  $\gamma_\omega(n)\le |S|^n$ for all $n$, so $\gamma_G\precsim e^n$.
\end{proof}

If there is a $d\in\N$ such that $\gamma_G\precsim n^d$, the group $G$
is of \emdef{polynomial growth} of degree at most $d$; if
$\gamma_G\sim e^n$, then $G$ is of \emdef{exponential growth};
otherwise $G$ is of \emdef{intermediate growth}. The existence of
groups of intermediate growth was first shown by
Grigorchuk~\cite{grigorchuk:growth}.


\begin{lem}\label{lem:sbgp}
  Let $H<G$ be an index-$N$ subgroup inclusion, let $\omega$ be a
  weight on $G$ with growth function $\gamma_\omega$, and let
  $\gamma^H_\omega$ denote the restriction of $\gamma_\omega$ to $H$:
  $$\gamma^H_\omega(n) = \#\{g\in H|\,\partial_\omega(g)\le n\}.$$
  Then there is a constant $K$ such that
  $$\gamma_\omega(n-K)\le N\gamma^H_\omega(n)\le\gamma_\omega(n+K)$$
  holds for all $n$.
\end{lem}
\begin{proof}
  Let $T$ be a transversal of $H$ in $G$, and set $K=\max_{t\in
    T}\partial_\omega(t)$. For every $g\in G$ of weight at most $n-K$,
  there is a unique $t\in T$ with $gt\in H$, and
  $\partial_\omega(gt)\le n$. The map $\{G\to H,\,g\mapsto gt\}$ is
  $N$-to-$1$, so its restriction to the set of elements of weight at
  most $n-K$ is at most $N$-to-$1$, proving the first inequality.

  Considering all $h\in H$ of weight at most $n$ and all $t\in T$, we have
  $N\gamma^H_\omega(n)$ distinct elements $ht\in G$, with
  $\partial_\omega(ht)\le n+K$. This proves the second inequality.
\end{proof}

\section{The Grigorchuk $2$-group}
Let $\Sigma^*$ be the set of finite sequences over $\Sigma=\{0,1\}$.
For $x\in\Sigma$ set $\overline x=1-x$. Define recursively the
following length-preserving permutations of $\Sigma^*$:
\begin{alignat*}{3}
  &&a(x\sigma)&=\overline x\sigma;\\
  b(0\sigma)&=0a(\sigma),&&&b(1\sigma)&=1c(\sigma);\\
  c(0\sigma)&=0a(\sigma),&&&c(1\sigma)&=1d(\sigma);\\
  d(0\sigma)&=0\sigma,&&&d(1\sigma)&=1b(\sigma).
\end{alignat*}
Then $G$, the Grigorchuk
$2$-group~\cite{grigorchuk:burnside,grigorchuk:gdegree}, is the group
generated by $S=\{a,b,c,d\}$.  It is readily checked that these
generators are of order $2$ and that $V=\{1,b,c,d\}$ is a Klein group.

Let $H=V^G$ be the normal closure of $V$ in $G$. It is of index $2$ in
$G$ and preserves the first letter of sequences; i.e.\ $H\cdot
x\Sigma^*\subset x\Sigma^*$ for all $x\in\Sigma$. There is an
injective homomorphism $\psi:H\to G\times G$, written $g\mapsto
(g_0,g_1)$, defined by $0g_0(\sigma)=g(0\sigma)$ and
$1g_1(\sigma)=g(1\sigma)$. As $H=\langle b,c,d,b^a,c^a,d^a\rangle$, we
can write $\psi$ explicitly as
$$\psi:\begin{cases}b\mapsto(a,c),\quad b^a\mapsto(c,a)\\ c\mapsto(a,d),\quad
  c^a\mapsto(d,a)\\ d\mapsto(1,b),\quad d^a\mapsto(b,1).\end{cases}$$

Let $B$ be the normal closure of $\langle b\rangle$ in $G$. We shall
use the following facts, whose proof appears for instance
in~\cite[Section~VIII.C]{harpe:cgt}: $B$ is of index $8$ in $G$,
and $\psi(H)=(B\times B)\rtimes\langle (a,1),(d,1)\rangle$ is of index
$8$ in $G\times G$, transversal to the order-$8$ dihedral group
$\langle (a,d),(d,a)\rangle$.

\section{The Growth of $G$}
\begin{defn}
  A weight $\omega$ on $(G,S)$ is \emdef{triangular} if for
  any ordering $(x,y,z)$ of $\{b,c,d\}$ one has
  $$\omega(x)\le\omega(y)+\omega(z).$$
\end{defn}
We shall only consider triangular weights on $G$, and implicitly use the
\begin{lem}\label{lem:nf}
  Let $\omega$ be a triangular weight.  Then every $g\in G$ admits a
  minimal form
  $$[*]a*a*a\dots*a[*],$$
  where $*\in\{b,c,d\}$ and the first and last $*$s are optional.
\end{lem}
\begin{proof}
  Let $w$ be a minimal form for $g\in G$. First, $w$ may not contain
  two equal consecutive letters, since they could be cancelled,
  shortening the representation of $g$. Second, two unequal
  consecutive letters among $\{b,c,d\}$ can be replaced by the third,
  and as $\omega$ is triangular this operation will shorten the length
  of $w$ while not enlarging its weight. Ultimately it will yield a
  word $w'$ of the required form, also representing $g$, and with
  smaller or equal weight. It can then be chosen as a minimal form
  of $g$.
\end{proof}

A lower estimate on the growth of $G$ comes from the following
proposition:
\begin{prop}\label{prop:cons}
  There are a weight $\omega$ on $S$ and constants $K\ge0,\eta<4$ such
  that for all $h\in H$, writing $\psi(h)=(h_0,h_1)$, we have
  \begin{equation}\label{eq:prop:cons}
    \partial_\omega(h) \le \eta\max\{\partial_\omega(h_0),\partial_\omega(h_1)\}+K.
  \end{equation}
\end{prop}

\begin{cor}\label{cor:main}
  The growth of $G$ satisfies
  $$\gamma_G\succsim e^{n^\alpha},$$
  where $\alpha = \log(2)/\log(\eta) \gneqq 0.5.$
\end{cor}
\begin{proof}
  Applying Lemma~\ref{lem:sbgp} (with constants $K_1$ and $K_2$) to
  the subgroup pairs $H<G$ and $\psi(H)<G\times G$, we obtain from the
  previous proposition
  $$\frac{\gamma_\omega(\eta x+K+K_1)}{2} \ge \gamma^H_\omega(\eta x+K)
  \ge \gamma^{\psi(H)}_{\omega\times\omega}(x) \ge \frac{\gamma_\omega(x-K_2)^2}{8}$$
  for all $x\ge K_2$.  At the cost of increasing $K$, we rewrite it as
  $\gamma_\omega(\eta x+K) \ge \gamma_\omega(x)^2/4$, whence by iteration
  $$\gamma_\omega\left(\eta^mx+\frac{\eta^m-1}{\eta-1}K\right)\ge\frac{\gamma_\omega(x)^{2^m}}{4^{2^m-1}}$$
  for all integers $m\ge1$ and $x\in\R_+$. Let $L\in\R_+$ be such that
  $\gamma(L)>4$. Now given $n\in\R_+$, let $m\in\N$ be maximal such
  that $x:=\eta^{-m}(n-K(\eta^m-1)/(\eta-1))$ be at least equal to $L$.
  We then have
  $\gamma_\omega(n)\ge(\gamma_\omega(x)/4)^{2^m}\ge(\gamma_\omega(L)/4)^{2^m}$,
  and since $m\approx\log(n)/\log(\eta)$, the corollary follows.
\end{proof}

Note that it is easy to prove Proposition~\ref{prop:cons} with the
constant $\eta=4$. Then by Corollary~\ref{cor:main} we would obtain
$\gamma_G\succsim e^{\sqrt n}$.
\begin{proof}[Proof of Proposition~\ref{prop:cons} for $\eta=4$]
  Consider the weight $\omega(s)=1$ for all $s\in S$ giving
  $\partial_\omega(g)=|g|$, and consider the following two
  homomorphisms (that they are homomorphisms was proven by Igor
  Lysionok~\cite{lysionok:pres}):
  \begin{align*}
    \sigma:&\begin{cases}G\to H\\
      a\mapsto c^a,\quad b\mapsto d,\quad d\mapsto c,\quad c\mapsto
      b;\end{cases}\\
    \tau:&\begin{cases}G\to G\\
      a\mapsto d,\quad b\mapsto 1,\quad d\mapsto a,\quad c\mapsto
      a.\end{cases}
  \end{align*}
  Then for any $g\in G$ we have $\psi(\sigma(g))=(\tau(g),g)$, and
  $\tau(g)$ is of length at most $4$, because $\tau(G)=\langle
  a,d\rangle$ is a dihedral group of order $8$, and
  $|\sigma(g)|\le2|g|+1$. Now for any $h\in H$, with
  $\psi(h)=(h_0,h_1)$, we have
  $$h=a\sigma(h_0)a\sigma(\tau(h_0)^{-1}h_1);$$
  indeed applying $\psi$ to the right-hand side gives
  $$(h_0,\tau(h_0))\big(\tau(\tau(h_0)^{-1}h_1),\tau(h_0)^{-1}h_1\big)=(h_0\tau(*),h_1)=(h_0,h_1),$$
  because $(\tau(G),1)\cap\psi(H)=(1,1)$.  Therefore we have
  \begin{align*}
    |h|&\le|a|+|\sigma(h_0)|+|a|+|\sigma\tau(*)|+|\sigma(h_1)|\\
    &\le 1+(2|h_0|+1)+1+8+(2|h_1|+1)\le4\max\{|h_0|,|h_1|\}+12.
  \end{align*}
\end{proof}

\section{Finite State Automata}
We prove Proposition~\ref{prop:cons} by constructing a ``finite state
automaton'' (see~\cite{hopcroft-u:automata} for an introduction) that,
given a pair $(h_0,h_1)$ of words in $\psi(H)$, constructs a preimage
$h$ satisfying~(\ref{eq:prop:cons}).

The automaton operates as follows: first each of the input words is
extended in an infinite word by ``padding'': extending it to the right
by a padding symbol $\pad$. At each step, the machine either outputs a
letter, or reads and deletes the first letter on each of the input
words. The concatenation of the output letters is the word
``computed'' by the machine.

Such a machine is most conveniently described as a directed,
labelled graph $\Gamma$. The vertices (called \emph{states}) of the
graph correspond to internal states of the machine, and the edges
(called \emph{transitions}) to input or output. Vertices are of two
types: \emph{input} and \emph{output}. The edges leaving an input
state are called \emph{input transitions} and are labelled by a pair
$(x,y)$ of input letters. The edges leaving an output state are called
\emph{output transitions} and are labelled by an output letter. There
are two distinguished vertices, the \emph{initial} and \emph{final}
states $\ast$ and $\dagger$.

Given a path $e=e_1\dots e_n$ in $\Gamma$, we write
$i(e)=(i(e)_0,i(e)_1)$ for the pair of words obtained by reading in
sequence along $e$ the labels of input transitions, and $o(e)$ for the
word obtained by reading the labels of output transitions.

Let $(u_0,u_1)$ be a pair of input words. Let $e$ be a path in
$\Gamma$ from $\ast$ to $\dagger$. If $i(e)$ is of the form
$(u_0\pad^*,u_1\pad^*)$, then $o(e)$ is a possible output of the
automaton.

Note that the output may not be uniquely determined --- if there is
always at most one output, the automaton is called
\emph{deterministic}. Note also that in some cases the automaton fails
to produce an output (if there exists no appropriate path). In our
situation these problems will not occur.

We now add an important condition: that the graph $\Gamma$ be finite.
(This means that the automaton has a bounded amount of available
memory.)
\begin{lem}\label{lem:autom}
  Let $\Gamma$ be a finite state automaton, and let $\partial$ be a
  weight on (input and output) words. Assume that for every directed
  loop $\ell$ in $\Gamma$ the following holds: $\partial o(\ell) <
  \eta/2(\partial i(\ell)_0+\partial i(\ell)_1)$.

  Then there is a constant $K$ such that, for every input pair
  $(u_0,u_1)$, there is an output word $v$ with
  $$\partial v \le \eta\max\{\partial u_0,\partial u_1\} + K.$$
\end{lem}
\begin{proof}
  Let $e$ be a path in $\Gamma$ from $\ast$ to $\dagger$ with
  $i(e)=(u_0\pad^*,u_1\pad^*)$. Since the graph is finite, say with
  $N$ vertices, $e$ may be decomposed as a product of paths
  $f_0\ell_1f_1\dots\ell_mf_m$, such that the $\ell_i$ are all loops
  in $\Gamma$ and $f_0f_1\dots f_m$ is of length at most $N$. Let $K$
  be the maximal weight of a word of length at most $N$. Then
  \begin{align*}
    \partial v &= \sum_{i=1}^m\partial o(\ell_i) + \sum_{i=0}^m\partial o(f_i)\\
    &\le \sum_{i=1}^m\eta/2(\partial i(\ell)_0+\partial i(\ell)_1) +
    K\\
    &\le \eta/2(\partial u_0+\partial u_1)+K \le \eta\max\{\partial u_0,\partial u_1\} + K.
  \end{align*}
\end{proof}

\section{Description of $\Gamma$}
We now describe the graph $\Gamma$ that computes $\psi$-preimages for
the proof of Proposition~\ref{prop:cons}. First we fix a minimal form
$\overline g$ for every $g\in G$, and denote by $\mf\subset S^*$ the
set of minimal forms. In Section~\ref{sec:construction} we will
describe a computer-aided process that constructs a triangular weight
$\omega:\{a,b,c,d\}\to\R_+^*$ and the finite oriented graph $\Gamma$.

Let us for now highlight the important properties of $\Gamma$:
\begin{itemize}
\item $V(\Gamma)\subset \mf\times\mf$, i.e.\ each vertex is
  identified with a pair of words $(u_0,u_1)$ in minimal form. Using
  this identification, $V(\Gamma)$ is also viewed as a subset of
  $G\times G$.
  
  The vertex set of $\Gamma$ corresponds to the machine's ``memory
  buffer'': roughly, when it reads a symbol pair, it adds it at the
  right of its buffer, while when it outputs a symbol, it deletes some
  part on the left of its buffer.
\item $\ast=\dagger=(\lambda,\lambda)$ is a vertex, where $\lambda$ is
  the empty word.
\item The definition in the previous section is slightly extended in
  that words (and not just letters) are allowed as edge-labels. Input
  transitions shall have labels of the form $(xa,ya)$ for all
  $x,y\in\{b,c,d\}$. Output vertices shall have a unique outbound
  edge, with a label in $\mf$.
\item An input transition labelled $(xa,ya)$ at a vertex $(u_0,u_1)$
  shall end at $(\overline{u_0xa},\overline{u_1ya})$.
  An output transition labelled $v$ at a vertex $(u_0,u_1)$ shall end at
  $(\overline{v_0u_0},\overline{v_1u_1})$.
\item At every input vertex $(u_0,u_1)\in\psi(H)$, and for every
  $u\in\mf\cap B$ of length of most $8$, there is an
  input transition labelled $(\pad^{|u|},u)$ from $(u_0,u_1)$ to a new
  vertex, directly followed by an output transition to
  $(\lambda,\lambda)$ labelled by some $v$ with $\psi(v)=(u_0,u_1u)$.
\item Except for these extra transitions, $\Gamma$ satisfies
  Lemma~\ref{lem:autom} with the same constant $\eta$.
\end{itemize}
For now, we suppose such a graph exists, and prove the proposition
under that assumption.

\begin{proof}[Proof of Proposition~\ref{prop:cons}]
  Let $(u_0,u_1)\in\psi(H)$ be a pair of words in minimal form; we are
  to construct a word $v$ such that $\psi(v)=(u_0,u_1)$
  satisfies~(\ref{eq:prop:cons}).

  First we may suppose that neither $u_0$ nor $u_1$ starts by `$a$';
  indeed we may add $b$'s at their beginning (remembering that $(b,1)$
  and $(1,b)$ are in $\psi(H)$) and construct a word $v'$ such
  that $\psi(v')=(b^iu_0,b^ju_1)$ with $i,j\in\{0,1\}$. The word
  $v=(ada)^id^jv'$ then satisfies $\psi(v)=(u_0,u_1)$ and is of weight
  at most $4$ more than $v'$, a fact that can be coped with by
  increasing the constant $K$.

  Second, we may suppose that $u_0$ is shorter than $u_1$, and that
  their lengths differ by at most $8$. Indeed the shorter of these two
  can be extended by copies of $(ad)^4$ which is trivial in $G$, and
  $u_1$ can be extended by an extra copy of $(ad)^4$. Again the cost
  of this operation is at worst an increase of $K$ by $8$.

  Now $\Gamma$ satisfies Lemma~\ref{lem:autom} for the constant
  $\eta$. This precisely means that there is an output word $v$
  satisfying~(\ref{eq:prop:cons}).
\end{proof}

\section{Construction of $\Gamma$}\label{sec:construction}
In this section we explain how a computer may be programmed to
construct the graph $\Gamma$ described in the previous section. The
construction itself involved a large amount of experimenting, to find
adequate values for the constants described below.

First an arbitrary triangular weight is selected, for instance
$\omega(s)\equiv1$ (in practice, values between $0.9$ and $1.1$ work
well); a tiny $\delta$ (about $0.01$), a largish $\eta'$ (around
$4$) and a limit-length $N$ (around $20$) are also chosen.

The \emph{quality} of an output transition $e$ labelled $v$ from
$(u_0,u_1)$ to $(u'_0,u'_1)$ is defined as
$$q(e)=\frac{\partial(u_0)+\partial(u_1)-\partial(u'_0)-\partial(u'_1)}{\partial v}
+\delta\left|\partial(u_0)-\partial(u_1)\right|-\delta\left|\partial(u'_0)-\partial(u'_1)\right|.$$
Therefore, edges connecting a vertex of large total weight to a vertex of
small total weight have a high quality. Edges connecting a vertex of
dissimilar weight to a vertex of similar-weight words also have high
quality.

The graph $\Gamma$ is now constructed iteratively. At all steps of the
iteration we shall have a graph $\Gamma$ satisfying all the required
conditions, except that it will have ``hanging edges'', that is, edges
not connected to any vertex. The purpose of the iteration process will
then be to attach the hanging edges, either to existent vertices in
$\Gamma$ or to new ones.

The graph $\Gamma$ starts with a single vertex, $(\lambda,\lambda)$,
and a hanging edge. Then, while there are hanging edges, the following
is performed:
\begin{itemize}
\item At each node $(u_0,u_1)$ with a hanging edge, all words
  $v\in\mf\cap H$ of length at most $N$ are tried; we write
  $\psi(v)=(v_0,v_1)$.
\item If the quality of the contracting edge from $(u_0,u_1)$ to
  $(\overline{v_0u_0},\overline{v_1u_1})$ is at least $1/\eta'$, then
  the type of $(u_0,u_1)$ is set to ``output'', and the hanging edge
  is replaced by an output transition from $(u_0,u_1)$ to
  $(\overline{v_0u_0},\overline{v_1u_1})$, labelled $v$.
\item If there is no such output transition of sufficient quality,
  then the type of $(u_0,u_1)$ is set to ``input'', and the hanging
  edge is replaced by $9$ edges to all the
  $(\overline{u_0xa},\overline{u_1ya})$ for all $x,y\in\{b,c,d\}$.
\end{itemize}

If the process stops, we have obtained a graph $\Gamma$ satisfying all
constraints, including Lemma~\ref{lem:autom} for the constant
$\eta'$. In fact, it may well satisfy Lemma~\ref{lem:autom} for
some smaller value of $\eta$: the triangular weight may be varied, and
it may happen that all edges of low quality are surrounded by edges of
high quality.

For this purpose, a second program was written. It takes as input the
description of $\Gamma$, and adjusts slightly the triangular
weight $\omega$; then it searches $\Gamma$ for loops, and on each loop
computes the weight of the input- and output-labels. If the adjustment
decreases the maximal ratio of input-weight to output-weight, it is
kept. Then another adjustment is tried, etc.

Finally the ``special transitions'' $(\pad^{|u|},u)$ are added to
$\Gamma$.

In actual computations, a graph $\Gamma$ with $160$ vertices (of which
$12$ are input states) was constructed. An appropriate weight was then
chosen to be
$$\omega(a)=1,\quad\omega(b)=3.33,\quad\omega(c)=2.8,\quad\omega(d)=1.06.$$
The resulting $\eta$ was $3.83414$.

\appendix
\section{The Graph $\Gamma$}
In this printout, the graph $\Gamma$ is described. The data are
structured in ``paragraphs'', each paragraph corresponding to a state
of type ``input''. Its nine successors are listed in the order
$(da,da),(da,ca),\dots,(ba,ba)$, and each of these successors is
identified by its type: if it is of type ``output'', it is followed by
the symbol
\begin{quote}
  \ttfamily OUTPUT($w$) $\rightarrow$ ($u_0$,$u_1$)
\end{quote}
meaning the word $w$ is output and the next state to be in is
$(u_0,u_1)$. If the successor is of type ``input'', it is followed by
the symbol \texttt{INPUT} and its successors are described in another
paragraph.

Input states are listed in lexicographical order, starting with pairs
$(u_0,u_1)$ of same length, then pairs with length-difference $1$,
then $2$, etc. In the lexicographical order $d$ comes before $c$ which
comes before $b$. The input states are $(\lambda,\lambda)$, $(da,da)$,
$(da,ca)$, $(ca,da)$, $(ca,ca)$, $(a,\lambda)$, $(da,d)$,
$(da,\lambda)$, $(ca,\lambda)$, $(ad,\lambda)$, $(ada,\lambda)$ and
$(aca,\lambda)$.

\newpage
{\footnotesize\ttfamily
\def\0#1{\parbox[t]{.75\linewidth}{\hrulefill\\ \rmfamily\textit{#1}}\\[1mm]}
\begin{tabbing}
\0{first, the initial (and terminal) vertex:}
($\lambda$,$\lambda$): INPUT $\rightarrow$ \=(da,da): INPUT\\
\>(da,ca): INPUT\\
\>(da,ba): OUTPUT(acad) $\rightarrow$ (a,$\lambda$)\\
\>(ca,da): INPUT\\
\>(ca,ca): INPUT\\
\>(ca,ba): OUTPUT(abad) $\rightarrow$ (a,$\lambda$)\\
\>(ba,da): OUTPUT(cada) $\rightarrow$ (a,$\lambda$)\\
\>(ba,ca): OUTPUT(bada) $\rightarrow$ (a,$\lambda$)\\
\>(ba,ba): OUTPUT(abad) $\rightarrow$ (da,$\lambda$)\\
\0{this finishes the description of $(\lambda,\lambda)$. Now comes the 
  description of those successors that are of type ``input''}
(da,da): INPUT $\rightarrow$ \=(dada,dada)=(adad,adad): OUTPUT(cacacaca) $\rightarrow$ ($\lambda$,$\lambda$)\\
\0{this is the first loop reached by the algorithm: its input-label is
  $(dada,dada)$ and its output-label is $cc^acc^a$. It has an
  input-output ratio of $\eta=3.69$. Note the simplification $dada=adad$}
\>(dada,daca)=(adad,daca): OUTPUT(cacabaca) $\rightarrow$ ($\lambda$,$\lambda$)\\
\>(dada,daba)=(adad,daba): OUTPUT(cacab) $\rightarrow$ (da,d)\\
\>(daca,dada)=(daca,adad): OUTPUT(acacabac) $\rightarrow$ ($\lambda$,$\lambda$)\\
\>(daca,daca): OUTPUT(acabacadadab) $\rightarrow$ ($\lambda$,$\lambda$)\\
\>(daca,daba): OUTPUT(cabacaca) $\rightarrow$ (aca,$\lambda$)\\
\>(daba,dada)=(daba,adad): OUTPUT(acacaba) $\rightarrow$ (da,d)\\
\>(daba,daca): OUTPUT(acabacac) $\rightarrow$ (aca,$\lambda$)\\
\>(daba,daba): OUTPUT(acabacacabacaca) $\rightarrow$ (ad,$\lambda$)\\
(da,ca): INPUT $\rightarrow$ \=(dada,cada)=(adad,cada): OUTPUT(bacacaca) $\rightarrow$ ($\lambda$,$\lambda$)\\
\>(dada,caca)=(adad,caca): OUTPUT(bacabaca) $\rightarrow$ ($\lambda$,$\lambda$)\\
\>(dada,caba)=(adad,caba): OUTPUT(bacab) $\rightarrow$ (da,d)\\
\>(daca,cada): OUTPUT(cabacacad) $\rightarrow$ ($\lambda$,$\lambda$)\\
\0{here is another loop at $(\lambda,\lambda)$: its input-label
  is $(daca,cada)$ and its output-label is $cb^acc^ad$. It has an
  output-input ratio of $\eta=3.04$}
\>(daca,caca): OUTPUT(acababab) $\rightarrow$ ($\lambda$,$\lambda$)\\
\>(daca,caba): OUTPUT(cabacacad) $\rightarrow$ (aca,$\lambda$)\\
\>(daba,cada): OUTPUT(cacadaba) $\rightarrow$ (da,d)\\
\>(daba,caca): OUTPUT(acabacab) $\rightarrow$ (aca,$\lambda$)\\
\>(daba,caba): OUTPUT(cabacacabacab) $\rightarrow$ (ad,$\lambda$)\\
(ca,da): INPUT $\rightarrow$ \=(cada,dada)=(cada,adad): OUTPUT(abacacac) $\rightarrow$ ($\lambda$,$\lambda$)\\
\>(cada,daca): OUTPUT(acabacacada) $\rightarrow$ ($\lambda$,$\lambda$)\\
\>(cada,daba): OUTPUT(acacadab) $\rightarrow$ (da,d)\\
\>(caca,dada)=(caca,adad): OUTPUT(abacabac) $\rightarrow$ ($\lambda$,$\lambda$)\\
\>(caca,daca): OUTPUT(cabacababadab) $\rightarrow$ ($\lambda$,$\lambda$)\\
\>(caca,daba): OUTPUT(cabacaba) $\rightarrow$ (aca,$\lambda$)\\
\>(caba,dada)=(caba,adad): OUTPUT(abacaba) $\rightarrow$ (da,d)\\
\>(caba,daca): OUTPUT(acabacacada) $\rightarrow$ (aca,$\lambda$)\\
\>(caba,daba): OUTPUT(acabacacabacaba) $\rightarrow$ (ad,$\lambda$)\\
(ca,ca): INPUT $\rightarrow$ \=(cada,cada): OUTPUT(acacacabada) $\rightarrow$ ($\lambda$,$\lambda$)\\
\>(cada,caca): OUTPUT(acabacabada) $\rightarrow$ ($\lambda$,$\lambda$)\\
\>(cada,caba): OUTPUT(acabadab) $\rightarrow$ (da,d)\\
\>(caca,cada): OUTPUT(cabacabad) $\rightarrow$ ($\lambda$,$\lambda$)\\
\>(caca,caca): OUTPUT(acabacabadabadab) $\rightarrow$ ($\lambda$,$\lambda$)\\
\>(caca,caba): OUTPUT(cabacabad) $\rightarrow$ (aca,$\lambda$)\\
\>(caba,cada): OUTPUT(cabadaba) $\rightarrow$ (da,d)\\
\>(caba,caca): OUTPUT(acabacabada) $\rightarrow$ (aca,$\lambda$)\\
\>(caba,caba): OUTPUT(badabababadababac) $\rightarrow$ (ad,$\lambda$)\\
\0{now come vertices whose lengths differ by $1$:}
(a,$\lambda$): INPUT $\rightarrow$ \=(ada,da): OUTPUT(caca) $\rightarrow$ (a,$\lambda$)\\
\>(ada,ca): OUTPUT(baca) $\rightarrow$ (a,$\lambda$)\\
\>(ada,ba): OUTPUT(b) $\rightarrow$ (da,da)\\
\>(aca,da): OUTPUT(caba) $\rightarrow$ (a,$\lambda$)\\
\>(aca,ca): OUTPUT(baba) $\rightarrow$ (a,$\lambda$)\\
\>(aca,ba): OUTPUT(b) $\rightarrow$ (ca,da)\\
\>(aba,da): OUTPUT(caba) $\rightarrow$ (da,$\lambda$)\\
\>(aba,ca): OUTPUT(baba) $\rightarrow$ (da,$\lambda$)\\
\>(aba,ba): OUTPUT(badac) $\rightarrow$ (a,$\lambda$)\\
(da,d): INPUT $\rightarrow$ \=(dada,dda)=(adad,a): OUTPUT(aca) $\rightarrow$ (ada,$\lambda$)\\
\>(dada,dca)=(adad,ba): OUTPUT(acad) $\rightarrow$ (ada,$\lambda$)\\
\>(dada,dba)=(adad,ca): OUTPUT(baca) $\rightarrow$ (ad,$\lambda$)\\
\>(daca,dda)=(daca,a): OUTPUT(aca) $\rightarrow$ (aca,$\lambda$)\\
\>(daca,dca)=(daca,ba): OUTPUT(acad) $\rightarrow$ (aca,$\lambda$)\\
\>(daca,dba)=(daca,ca): OUTPUT(cabacacad) $\rightarrow$ (ad,$\lambda$)\\
\>(daba,dda)=(daba,a): OUTPUT(acabadab) $\rightarrow$ ($\lambda$,$\lambda$)\\
\>(daba,dca)=(daba,ba): OUTPUT(acacadab) $\rightarrow$ ($\lambda$,$\lambda$)\\
\>(daba,dba)=(daba,ca): OUTPUT(badabaca) $\rightarrow$ (aca,$\lambda$)\\
\0{now come vertices whose lengths differ by $2$:}
(da,$\lambda$): INPUT $\rightarrow$ \=(dada,da)=(adad,da): OUTPUT(caca) $\rightarrow$ (ad,$\lambda$)\\
\>(dada,ca)=(adad,ca): OUTPUT(baca) $\rightarrow$ (ad,$\lambda$)\\
\>(dada,ba)=(adad,ba): OUTPUT(acad) $\rightarrow$ (ada,$\lambda$)\\
\>(daca,da): OUTPUT(cabacaca) $\rightarrow$ (ad,$\lambda$)\\
\>(daca,ca): OUTPUT(cabacacad) $\rightarrow$ (ad,$\lambda$)\\
\>(daca,ba): OUTPUT(acad) $\rightarrow$ (aca,$\lambda$)\\
\>(daba,da): OUTPUT(badabaca) $\rightarrow$ (ada,$\lambda$)\\
\>(daba,ca): OUTPUT(badabaca) $\rightarrow$ (aca,$\lambda$)\\
\>(daba,ba): OUTPUT(acacadab) $\rightarrow$ ($\lambda$,$\lambda$)\\
(ca,$\lambda$): INPUT $\rightarrow$ \=(cada,da): OUTPUT(acacada) $\rightarrow$ (ad,$\lambda$)\\
\>(cada,ca): OUTPUT(acabada) $\rightarrow$ (ad,$\lambda$)\\
\>(cada,ba): OUTPUT(abad) $\rightarrow$ (ada,$\lambda$)\\
\>(caca,da): OUTPUT(cabacaba) $\rightarrow$ (ad,$\lambda$)\\
\>(caca,ca): OUTPUT(cabacabad) $\rightarrow$ (ad,$\lambda$)\\
\>(caca,ba): OUTPUT(abad) $\rightarrow$ (aca,$\lambda$)\\
\>(caba,da): OUTPUT(badababa) $\rightarrow$ (ada,$\lambda$)\\
\>(caba,ca): OUTPUT(badababa) $\rightarrow$ (aca,$\lambda$)\\
\>(caba,ba): OUTPUT(abacadab) $\rightarrow$ ($\lambda$,$\lambda$)\\
(ad,$\lambda$): INPUT $\rightarrow$ \=(adda,da)=($\lambda$,da): INPUT\\
\0{in fact, the destination edge is $(da,\lambda)$, since we list
  edge pairs with longest word first. This is irrelevant, since word
  reductions are similar when operating on $(u_0,u_1)$ and ond $(u_1,u_0)$}
\>(adda,ca)=($\lambda$,ca): INPUT\\
\>(adda,ba)=($\lambda$,ba): OUTPUT(ada) $\rightarrow$ (a,$\lambda$)\\
\>(adca,da)=(aba,da): OUTPUT(caba) $\rightarrow$ (da,$\lambda$)\\
\>(adca,ca)=(aba,ca): OUTPUT(baba) $\rightarrow$ (da,$\lambda$)\\
\>(adca,ba)=(aba,ba): OUTPUT(badac) $\rightarrow$ (a,$\lambda$)\\
\>(adba,da)=(aca,da): OUTPUT(caba) $\rightarrow$ (a,$\lambda$)\\
\>(adba,ca)=(aca,ca): OUTPUT(baba) $\rightarrow$ (a,$\lambda$\\
\>(adba,ba)=(aca,ba): OUTPUT(b) $\rightarrow$ (ca,da)\\
\0{finally some vertices whose lengths differ by $3$:}
(ada,$\lambda$): INPUT $\rightarrow$ \=(adada,da)=(dad,da): OUTPUT(acac) $\rightarrow$ (ada,$\lambda$)\\
\>(adada,ca)=(dad,ca): OUTPUT(acab) $\rightarrow$ (ada,$\lambda$)\\
\>(adada,ba)=(dad,ba): OUTPUT(acad) $\rightarrow$ (ad,$\lambda$)\\
\>(adaca,da): OUTPUT(c) $\rightarrow$ (daca,a)\\
\>(adaca,ca): OUTPUT(b) $\rightarrow$ (daca,a)\\
\>(adaca,ba): OUTPUT(b) $\rightarrow$ (daca,da)\\
\>(adaba,da): OUTPUT(c) $\rightarrow$ (daba,a)\\
\>(adaba,ca): OUTPUT(b) $\rightarrow$ (daba,a)\\
\>(adaba,ba): OUTPUT(b) $\rightarrow$ (daba,da)\\
(aca,$\lambda$): INPUT $\rightarrow$ \=(acada,da): OUTPUT(caba) $\rightarrow$ (ada,$\lambda$)\\
\>(acada,ca): OUTPUT(baba) $\rightarrow$ (ada,$\lambda$)\\
\>(acada,ba): OUTPUT(acacadab) $\rightarrow$ (ad,$\lambda$)\\
\>(acaca,da): OUTPUT(caba) $\rightarrow$ (aca,$\lambda$)\\
\>(acaca,ca): OUTPUT(baba) $\rightarrow$ (aca,$\lambda$)\\
\>(acaca,ba): OUTPUT(b) $\rightarrow$ (caca,da)\\
\>(acaba,da): OUTPUT(cababadab) $\rightarrow$ ($\lambda$,$\lambda$)\\
\>(acaba,ca): OUTPUT(bababadab) $\rightarrow$ ($\lambda$,$\lambda$)\\
\>(acaba,ba): OUTPUT(b) $\rightarrow$ (caba,da)\\
\end{tabbing}}

\def\nop#1{}\font\cyr=wncyr8
\providecommand{\bysame}{\leavevmode\hbox to3em{\hrulefill}\thinspace}

\end{document}